\documentclass[10pt]{amsart}
\usepackage{amsmath,amssymb}
\let\newpf\proof \let\proof\relax 
\newenvironment{pf}{\newpf[\proofname]}{\qed\endtrivlist}

\def\be{\begin{equation}}
\def\ee{\end{equation}}

\def\ba{{\begin{align}}}
\def\ea{{\end{align}}}

\def\u{{\mathbb U}}

\def\bm{\begin{matrix}}
\def\em{\end{matrix}}

\def\bp{\begin{pmatrix}}
\def\ep{\end{pmatrix}}

\def\SL{\mathrm{SL}}
\def\PSL{\mathrm{PSL}}
\def\SO{\mathrm{SO}}

\def\0{{\mathbf 0}}

\newtheorem{thm}{Theorem}[section]
\newtheorem{cor}[thm]{Corollary}

\newtheorem{lemma}[thm]{Lemma}

\theoremstyle{remark}
\newtheorem{rem}{Remark}[section]

\numberwithin{equation}{section}

\def \bn {\hfill \\ \smallskip\noindent}

\theoremstyle{definition}

\def\proof{\bn {\bf Proof.} }

\def\ssm{\smallsetminus}
\def\tr{{\text{tr}}}
\renewcommand{\setminus}{\ssm}

\newcommand{\id}{\operatorname{id}}

\newcommand{\C}{{\mathbb C}}
\newcommand{\D}{{\mathbb D}}

\newcommand{\Q}{{\mathbb Q}}
\newcommand{\R}{{\mathbb R}}

\newcommand{\Z}{{\mathbb Z}}

\renewcommand{\i}{{\bar i}}

\def\B0{{\bold{0}}}

\catcode`\@=12

\def\Empty{}
\newcommand\oplabel[1]{
  \def\OpArg{#1} \ifx \OpArg\Empty {} \else
  	\label{#1}
  \fi}
		
\newcommand{\comm}[1]{}
\newcommand{\comment}[1]{}

\begin{document}

\bigskip\bigskip

\title[Almost reducibility and absolute continuity I]
{Almost reducibility and absolute continuity I}
\author {Artur Avila}

\address{
CNRS UMR 7586, Institut de Math\'ematiques de Jussieu\\
175 rue du Chevaleret\\
75013, Paris, France \&
IMPA, Estrada Dona Castorina 110\\
22460-320, Rio de Janeiro, Brazil.
}
\email{artur@math.jussieu.fr}

\date{\today}

\begin{abstract}

We consider one-frequency analytic $\SL(2,\R)$ cocycles.  Our main result
establishes the Almost Reducibility Conjecture in the case of exponentially
Liouville frequencies.  Together with our earlier work, this implies that
all cocycles close to constant are almost reducible, independent of the
frequency.  In our forthcoming work, we discuss applications to the analysis
of the absolutely continuous spectrum of one-frequency Schr\"odinger
operators.

\end{abstract}

\comm{
We consider one-dimensional quasiperiodic Schr\"odinger operators with
$(Hu)_n=u_{n+1}+u_{n-1}+v(\theta+n\alpha) u_n$ where $v:\R/\Z \to \R$ is
analytic and $\theta,\alpha \in \R$.  We prove that if $v$ is sufficiently
small (independent of $\alpha$ and $\theta$) then the spectral measures are
absolutely continuous.  In view of previous work (Avila-Jitomirskaya,
Avila), we concentrate here on the case where $\alpha$ is very well
approximated by rational numbers.  The proof is based on a very precise
control of the dynamics of the associated cocycle (``almost reducibility'').
}

\setcounter{tocdepth}{1}

\maketitle

\section{Introduction}

Here we consider one-frequency analytic $\SL(2,\R)$ cocycles, that is,
linear skew-products over an irrational rotation $x \mapsto x+\alpha$ of
the circle $\R/\Z$ which have the form $(\alpha,A):(x,w) \mapsto
(x+\alpha,A(x) \cdot w)$ with $A:\R/\Z \to \SL(2,\R)$.  The iterates of the
cocycle have the form $(\alpha,A)^n=(n \alpha,A_n)$ with $A_n(x)=A(x+(n-1)
\alpha) \cdots A(x)$, and the Lyapunov exponent is defined by
\be \label {L}
L=\lim_{n \to \infty} \frac {1} {n} \int \ln \|A_n(x)\| dx.
\ee
We say that $(\alpha,A)$ is uniformly hyperbolic if the cocycle iterates
grow exponentially uniformly on $x \in \R/\Z$.  Uniform hyperbolicity is
robust (it corresponds to an open set of cocycles) and easily analyzed.

Recently, an approach to the global theory of one-frequency cocycles has
been proposed \cite {A2}, \cite {A3}.
In it, cocycles which are not uniformly hyperbolic
are classified in three regimes:
\begin{enumerate}
\item Supercritical (or nonuniformly hyperbolic),
if the Lyapunov exponent is positive but $(\alpha,A)$ is not uniformly
hyperbolic,
\item Subcritical, if the cocycle iterates $\|A_n(z)\|$ are uniformly
subexponentially bounded through some strip $\{|\Im z|<\epsilon\}$,
\item Critical otherwise.
\end{enumerate}
A key point of this classification is that
(in the complement of uniform
hyperbolicity) both supercriticality and
subcriticality are stable (respectively, by \cite {BJ1} and \cite {A2}),
while criticality is unstable (it is the boundary
of supercriticality, see \cite {A3}).  Moreover, in
\cite {A3} it is shown that criticality
is ``negligible'' in the sense
that it does not appear at all in typical one-parameter families
(this is quite convenient for the theory since very little is known about
the dynamics of critical cocycles, apart that they are rare).
Naturally, one still is left with the problem of describing the stable
regimes.

As it turns out, the understanding of supercriticality is quite developed,
through the works of Bourgain, Goldstein, Jitomirskaya and Schlag
\cite {BG}, \cite {GS}, \cite {BJ1},
\cite {GS1}, \cite {GS2}.  Subcriticality on the
other hand, is a relatively new concept, which was first suggested to be
relevant in 2006 (see \cite {AJ}).
In fact, \cite {AJ} basically proposed that
the well developed theory of cocycles close to constant (\cite {E},
\cite {BJ2}, \cite {AJ}, \cite {AFK}) can be applied to
all subcritical cocycles, by the application of suitable coordinate changes.
Recall that in the cocycle context the natural notion of coordinate change
is given by a conjugacy $(x,w) \mapsto (x,B(x) \cdot w)$ with $B:\R/\Z \to
\PSL(2,\R)$ analytic, which takes $(\alpha,A)$ to $(\alpha,A')$ with
$A'(x)=B(x+\alpha) A(x) B(x)^{-1}$.  Let us
say that $(\alpha,A)$ is {\it almost reducible} if there exist $\epsilon>0$
and a sequence
of analytic maps $B^{(n)}:\R/\Z \to \PSL(2,\R)$, admitting holomorphic
extensions to the common strip $\{|\Im z|<\epsilon\}$ such that
$B^{(n)}(z+\alpha) A(z) B^{(n)}(z)^{-1}$ converges to a constant uniformly in
$\{|\Im z|<\epsilon\}$ (the $B^{(n)}$ themselves are allowed to diverge).
Essentially by definition, the concept of almost reducibility ``captures''
the domain of applicability of ``local theories'' of cocycles close to
constant.  The basic hope expressed by \cite {AJ} can be thus expressed in
the form of the {\it Almost Reducibility Conjecture} (ARC):
subcriticality implies almost reducibility.

Our first main result establishes a generic version of the ARC.
Let us say that $\alpha \in \R \setminus \Q$ is
exponentially Liouville if $\limsup \frac {\ln q_{n+1}} {q_n}>0$, where
$q_n$ is the sequence of denominators of continued fraction approximants.

\begin{thm} \label {thm1}

If $\alpha \in \R \setminus \Q$ be exponentially Liouville, then any
subcritical cocycle $(\alpha,A)$ is almost reducible.

\end{thm}

\begin{rem}

We should mention that the results of this paper (obtained in 2006-2007),
preced the results of \cite {A2} and \cite {A3} (obtained in 2008-2009),
and (being at the time the only evidence for the ARC),
played a large role in motivating those works.  On the other
hand, we have recently established the ARC
for almost every frequency (by very different methods, in particular
making use of \cite {A2}).

\end{rem}

Together with our previous analysis of frequencies $\alpha$ which are not
exponentially Liouville, Theorem \ref {thm1} implies:

\begin{cor} \label {cor2}

Any one-frequency cocycle close to constant is almost reducible.

\end{cor}

\begin{pf}

Any one-frequency cocycle which is close to constant is either
uniformly hyperbolic or subcritical (this is essentially due to \cite {BJ1}
and \cite {BJ2}, and it is explicitly obtained in \cite {A1} by a different
argument).  Uniformly hyperbolic
cocycles are always almost reducible: they can be conjugated to a diagonal
cocycle, which can then be conjugated arbitrarily close to a constant one
using approximate solutions of the cohomological equation.  If $\alpha$ is
exponentially Liouville, the result then follows by Theorem \ref {thm1}.

The complementary case was established earlier in \cite {A1}.\footnote
{For $\alpha$ Diophantine, i.e., under the condition $\ln q_{n+1}=O(\ln q_n)$,
this was established in \cite {AJ} Theorem 4.1 (rigorously speaking, \cite
{AJ} only deals with the case of Schr\"odinger
cocycles, but this case implies the general one by a simple abstract
argument, see Lemma 2.2 of \cite {AJ2}).  Under the weaker condition $\ln
q_{n+1}=o(q_n)$ (defining the complement of the exponentially Liouville
regime), almost reducibility near constants follows from \cite {A1}, which
provides the necessary estimates for the argument of \cite {AJ}
(see Section 3.7 of \cite {A1}).}
\end{pf}

\begin{cor}

Almost reducibility is stable, in the sense that it defines an open set in
$(\R \setminus \Q) \times C^\omega(\R/\Z,\SL(2,\R))$.

\end{cor}

\begin{pf}

Let $(\alpha,A)$ be almost reducible and let $B^{(n)}$ be the sequence of
conjugacies as in the definition.
Let $(\alpha^{(n)},A^{(n)})$ be any sequence of non-almost reducible
cocycles converging to $(\alpha,A)$.
Then there exists a sequence $j_n \to \infty$ such that
$\tilde A^{(n)}(x)=B^{(n)}(x+\alpha^{(j_n)}) A^{(j_n)}(x) B^{(n)}(x)^{-1}$
converges to a constant.  By Corollary \ref {cor2}, $\tilde A^{(n)}$
must be almost reducible for large $n$.
Since almost reducibility is conjugacy invariant,
$(\alpha^{(j_n)},A^{(j_n)})$ is almost reducible, contradiction.
\end{pf}

{\bf Acknowledgements:} This research was partially conducted during the
period the author was a Clay Research Fellow.

\section{Rational approximation}

Let $C^\omega_\delta(\R,*)$, be the space of bounded analytic functions
with values in $*=\R,\SL(2,\R),...$, which admit a bounded analytic
extension to the strip
$\{|\Im z|<\delta\}$, with the norm
$\|a\|_\epsilon=\sup_{|\Im z|<\delta} |a(z)|$.  We let $C^\omega(\R/\Z,*)
\subset C^\omega(\R,*)$ be the subspace of $1$-periodic functions, with the
same norm.  Let $R_\theta=\bp \cos 2 \pi \theta & -\sin 2 \pi \theta\\
\sin 2 \pi \theta&\cos 2 \pi \theta \ep$.

Theorem \ref {thm1} will be obtained as a consequence of an estimate for
{\it periodic} cocycles with large period.

\begin{thm} \label {q}

For every $0<\epsilon<\epsilon_0$ there exists $C>0$ such that if
$\delta_1>0$ is sufficiently small, then for every $p/q \in \Q$ with
$q$ sufficiently large, if
$A \in C^\omega_{\epsilon_0}(\R/\Z,\SL(2,\R))$ is
such that
\be \label {cond}
\max_{0 \leq k \leq q}
\ln \|A_k\|_{\epsilon_0} \leq \delta_1 q, \quad \text {with} \quad
A_k(x)=A(x+(k-1) p/q) \cdots A(x)
\ee
then there exists
$B \in C^\omega_\epsilon(\R/\Z,\PSL(2,\R))$ and $R_* \in \SO(2,\R)$
such that $\|B(x)\|_\epsilon \leq e^{C \delta_1 q}$, and
$\|B(z+p/q)A(z)B(z)^{-1}-R_*\|_\epsilon \leq e^{-\delta_1 q}$.

\end{thm}

\noindent{\it Proof of Theorem \ref {thm1}.}
Let $(\alpha,A)$ be
subcritical.  By definition, there exists $\epsilon_0>0$ such that
\be
\lim_{n \to \infty} \sup_{|\Im z|<\epsilon_0} \frac {1} {n} \ln
\|A(z+(n-1) \alpha) \cdot A(z)\|=0.
\ee
Thus for every $\delta>0$, there exists $n \geq 1$ such that
\be
\sup_{0 \leq k \leq n} \sup_{|\Im z|<\epsilon_0}
\|A(z+(k-1) \alpha) \cdots A(z)\| \leq e^{\delta n}.
\ee
In particular, if $q \geq n$ and
$p/q$ is close to $\alpha$, then
\be \label {cond1}
\sup_{0 \leq k \leq n} \|A_k\|_{\epsilon_0} \leq 2 \delta n,
\quad \text {with} \quad A_n(x)=A(x+(k-1) p/q) \cdots A(x),
\ee
which implies (\ref {cond}) with $\delta_1=2\delta$.

Assume now that $\alpha$ is exponentially Liouville.  Then there exists
$\delta'>0$ such that we may choose
$p/q$ arbitrarily close to $\alpha$ and satisfying
$|\alpha-\frac {p} {q}|<e^{-\delta' q}$.  Fix $0<\epsilon'<
\epsilon<\epsilon_0$ and
let $C$ be as in Theorem \ref {q}.  Select $0<\delta_1<\frac {1} {10 C}
\delta'$ and
letting $B$ be and $R_*$ be as in Theorem \ref {q}, we get
\begin{align}
\|B(z+\alpha) A(z) &B(z)^{-1}-R_*\|_{\epsilon'} \leq\\
\nonumber
& \|B(z+p/q) A(z) B(z)^{-1}-R_*\|_\epsilon+\|A\|_\epsilon
\|B\|_\epsilon |\alpha-\frac {p} {q}| \|\partial B\|_{\epsilon'}\\
\nonumber
& \leq \|B(z+p/q) A(z) B(z)^{-1}-R_*\|_\epsilon+\|A\|_\epsilon
\|B\|_\epsilon^2 C(\epsilon,\epsilon') |\alpha-\frac {p} {q}|\\
\nonumber
& \leq C(\epsilon,\epsilon') e^{-\delta' q/2}.
\end{align}
Since $q$ can be taken arbitrarily large, the result follows.
\qed

We can actually obtain much more information than what is
described above.  For further applications (see
\cite {A4}), we will need such stronger estimates, but only in
a particular case (which is not the hardest to prove).

\begin{thm} \label {qsimpler}

For every $0<\epsilon<\epsilon_0$ there exist $\delta_2>0$ and $C>0$
such that if
$\delta_1>0$ is sufficiently small, then for every $p/q \in \Q$ with
$q$ sufficiently large, if
$A \in C^\omega_{\epsilon_0}(\R/\Z,\SL(2,\R))$ is
is such that (\ref {cond}) holds and
$t=\tr A_q$ satisfies
$|\hat t_0|<2$ and $2-|\hat t_0| \geq e^{-\delta_1 q}$
then there exists
$B \in C^\omega_\epsilon(\R/\Z,\PSL(2,\R))$ and $\theta \in
C^\omega_\epsilon(\R/\Z,\R)$ such that $\theta$ is $1/q$-periodic,
$\|\theta-\hat \theta_0\|_\epsilon \leq e^{-\delta_2 q}$ and
$\|B(x)\|_\epsilon \leq e^{C \delta_1 q}$, and
$B(z+p/q)A(z)B(z)^{-1}=R_{\theta(z)}$.

\end{thm}

\comm{
In the setup of Theorem \ref {q}, there exists
$\delta_2=\delta_2(\epsilon_0,\epsilon,C)$ such that $\delta_1>0$ and $q_0$
let $t(x)=\tr A_q(x)$,
and assume that $|\hat t_0|<2$ and $-\ln (2-|\hat t_0|) \leq \delta_1 q$.
Then $B$ can be taken in $C^\omega_\epsilon(\R/\Z,\SL(2,\R))$ and
$\tilde A(x)=R_{\theta(x)}$ with $\theta \in C^\omega_\epsilon(\R/\Z,\R)$
which is $1/q$-periodic and satisfies
$\|\theta-\hat \theta_0\|_\epsilon<e^{-\delta_0 q}$, where
$\delta=\delta(\epsilon,\epsilon_0)$.

\end{thm}
}

\section{Proof of Theorem \ref {qsimpler}}

\subsection{Preliminary estimates}

\begin{lemma} \label {Pu=0}

Let $P(x)=\bp a(x)&b(x)\\c(x)&d(x) \ep$ with $a,b,c,d \in
C^\omega_{\epsilon_0}(\R/\Z,\C)$.  Assume that $\det
P$ is identically vanishing and
$\delta \leq \|P(z)\| \leq 1$ through $\{|\Im z|<\epsilon_0\}$.
Then there exists $u \in C^\omega_\epsilon(\R/\Z,\C^2)$
such that $P(z) u(z)=0$ and $C^{-1}
\delta^C \leq \|u(z)\| \leq 1$ through $\{|\Im z|<\epsilon\}$.
Here $C=C(\epsilon_0,\epsilon)$.

\end{lemma}

\begin{pf}

If $c$ or $d$ vanishes identically, the result is obvious.  Indeed, if $c$
vanishes identically, for instance, then either $a$ vanishes identically
(and $u=(1,0)$ will do) or $d$ vanishes identically (and $u=(-b,a)$ will
do).

Let us assume that both $c$ and $d$ are not identically vanishing.
Define a meromorphic function (not identically $\infty$)
$\phi(x)=\frac {a(x)} {b(x)}=\frac {c(x)} {d(x)}$.  All estimates below are
for $|\Im x|<\epsilon$, and $C=C(\epsilon_0,\epsilon)$.

If $1/4<|\phi(x)|<1$ then $|D \phi(x)| \leq C/\delta$.  Thus
the $C^{-1} \delta$-neighborhood of $\phi(\{|\Im x|=\epsilon\})$ intersects
$\{1/2<|\kappa|<3/4\}$ in a set of $\kappa$ of Lebesgue measure at most
$1/10$.
This implies that there exists $|\kappa|<3/4$ such that
$|\phi(x)-\kappa|>C^{-1} \delta$ for every $x$ with $|\Im x|=\epsilon$,
and such that for every $y$ with $|\Im y|<\epsilon$ and
$\phi(y)=\kappa$ we have $|D \phi(y)|>C^{-1}$.  Up to replacing $P$
by $P \left (\bm 1&-\overline \kappa\\ -\kappa & 1 \em \right )$, we may
suppose that $\kappa=0$.  In particular, the zeros of $\phi$ are simple. 
Let us estimate the number of zeros of $\phi$ in $|\Im x|<\epsilon$.

If $\phi(x)=0$, then either for $\psi_0=a$ and $\psi_1=c$ or $\psi_0=b$ and
$\psi_1=d$ we have $|D \psi_0(x)|>C^{-1} \delta$, $|\psi_1(x)|>C^{-1}
\delta$.  This implies that we can cover the zeros of $\phi$ in
$\{|\Im x|<\epsilon\}$ with disjoint disks $D$
of radius $C^{-1} \delta$, such that $\max_{\psi=a,b}
\inf_{x \in \partial D} |\psi(x)|>C^{-1} \delta^2$.  The zeros (of $a$ or
$b$) in such disks
persist truncation of the Fourier series keeping frequencies at most
$-C \ln \delta$, hence $\phi$ has at most
$-C \ln \delta$ zeros in $|\Im x|<\epsilon$.

Let $p(z)=\prod_{s=1}^N (z-z_s)$ where $z_s=e^{2 \pi i \theta_s}$ and
$\theta_s$, $1 \leq s \leq N$ are the zeros of $\phi$.
Let $u_1(\theta)=p(e^{2 \pi i \theta})$, and let
$u_2(x)=u_1(x)/\phi(x)$.  Since the zeros of $\phi$ are simple, $u_1$ and
$u_2$ are bounded holomorphic functions in $|\Im x|<\epsilon$.  Let
$\lambda=\left \| \bp u_1 \\ u_2 \ep \right \|_\epsilon$.
We claim that $u=\lambda^{-1}
\left (\bm -u_2\\u_1 \em \right )$ has the desired properties.

Clearly $-a u_2+b u_1=u_2 (-a+b \phi)=0$ and similarly $-c u_2+d u_1=0$, so
that $P u=0$.  We also have $\|u(x)\|=1$ in $|\Im x|<\epsilon$.  We need to
show that $C^{-1} \delta^C \leq \|(u_1(x),u_2(x))\| \leq C \delta^{-C}$ in
$|\Im x|<\epsilon$.

Since the number $N$
of zeros of $\phi$ in $|\Im x|<\epsilon$ is bounded by $-C
\ln \delta$, we have $|u_1(x)| \leq C \delta^{-C}$ in $|\Im
x|<\epsilon$, and since $|\phi(x)|>C^{-1} \delta$ in $|\Im x|=\epsilon$, we
also have $|u_2(x)| \leq C \delta^{-C}$ in $|\Im x|<\epsilon$.  This gives
the upper estimate.

To conclude, let us show that $a(x)/u_1(x) \leq C^{-1} \delta^C$ and
$c(x)/u_1(x) \leq C^{-1} \delta^C$ for $|\Im x|<\epsilon$.  This implies the
lower estimate, since $(a,b)$ and $(c,d)$ are multiples of $(u_1,u_2)$ and
$\|P(x)\| \geq \delta$.

Since $a(x)/u_1(x)$ and $c(x)/u_1(x)$
are holomorphic in $|\Im x|<\epsilon_0$ and $\|P(x)\| \leq 1$ in $|\Im
x|<\epsilon_0$, it is enough to show that $u_1(x) \geq C \delta^{-C}$ for
$|\Im x|=\epsilon_0$.  But clearly $|u_1(x)| \geq
|e^{\mp 2 \pi \epsilon}-e^{\mp 2 \pi \epsilon_0}|^N$ if $\pm \Im
x=\epsilon_0$, where $N<-C \ln \delta$ is the number of zeros of $\phi$ in
$|\Im x|<\epsilon$.  This implies the lower estimate.
\comm{
For every
$\epsilon<|t|<\epsilon_0$, $\inf_{\Im x=t} |u_1(x)| \geq \prod |e^{-2 \pi t}-
e^{-2 \pi t_s}|$, where $t_s=\Im \theta_s$.  The set $K$
of $z \in \C$ such that
$\prod |z-e^{-2 \pi t_s}|<C^{-1} \delta^C$ has analytic capacity less that
$C^{-1}$.  Thus there exists $\epsilon<t_\pm<\epsilon_0$ such that $e^{\mp 2
\pi t_\pm} \notin K$.
}
\comm{
If $P$ is real-symmetric, the above proof does not immediately gives a
real-symmetric $u$ because we needed to replace $P$ with
$P \left (\bm 1&-\overline \kappa\\ -\kappa & 1 \em \right )$ at one step,
which does not preserve real-symmetry.
To solve this problem, let us first replace it by
$\left (\bm -1 & i\\1 & i \em \right )
P \left (\bm 1&-1\\i&i \em \right )$.  Then the real-symmetry assumption
gets replaced by the conditions $d(x)=\overline {a(x)}$ and
$c(x)=\overline {b(x)}$ for $x \in \R$, and replacing $P$ by
$P \left (\bm 1&-\overline \kappa\\ -\kappa & 1 \em \right )$ does not
affect it.  The condition $\det P=0$ implies $|a(x)|=|b(x)|$ for $x \in \R$,
and so $|\phi(x)|=1$ for $x \in \R$.
}
\end{pf}

\begin{lemma} \label {mu}

For every $0<\epsilon<\epsilon_2$, there exists $\delta>0$ and $C>0$ such if
$\delta_1$ is sufficiently small and $q$ is sufficiently large, the
following property holds.
Let $\mu \in C^\omega_{\epsilon_2}(\R/\Z,\C)$,
and let $\mu_k=\prod_{j=0}^{k-1} \mu(x+j p/q)$.  Assume that
$\|\mu_k\|_{\epsilon_2} \leq e^{\delta_1 q}$,
$1 \leq k \leq q$ and that $\|\mu_q^{-1}\|_{\epsilon_2} \leq
e^{\delta_1 q}$.\footnote {An hypothesis such as $\|\mu_q^{-1}\|_{\epsilon'_2} \leq
e^{o(q)}$, $\epsilon'_2$ close to $\epsilon_2$, which is enough for our
purposes, follows from the assumption that $\sup_{x \in \R/\Z}
|\mu_q(x)-\lambda|<e^{-\kappa
q}$ for some $\lambda \in \C$ with $|\lambda| \geq e^{-o(q)}$, by convexity.}
Then there exist $\psi,\theta \in
C^\omega_\epsilon(\R/\Z,\C)$ such that
\be
\mu(z)=e^{2 \pi i \theta(z)}
\frac {e^{2 \pi i \psi(z+p/q)}} {e^{2 \pi i\psi(z)}}
\ee
and $\|\psi\|_\epsilon \leq C \delta_1 q$, $\ln \|\theta-\hat
\theta_0\|_\epsilon \leq -\delta q$ and $\theta$ is $1/q$-periodic.
Moreover,
\begin{enumerate}
\item If $|\mu(x)|=1$ for every $x \in \R$, then $\Im \psi(x)=0$ for every
$x \in \R$,
\item If $\mu(x) \in \R$ for every $x \in \R$, then $\Re \psi(x)=0$ for
every $x \in \R$.
\end{enumerate}

\end{lemma}

\begin{pf}

Fix $\epsilon<\epsilon_3<\epsilon_2'<\epsilon_2$.
Let $\mu(x)=e^{2 \pi i (d x+\phi(x))}$, where $\phi \in
C^\omega_{\epsilon_2'}(\R/\Z,\C)$.
Then $\mu_q(x)=e^{2 \pi i (d q x)+\tilde \phi(x)}$ with $\tilde \phi \in
C^\omega_{\epsilon_2'}(\R/\Z,\C)$.  Let $\lambda$ be the average of $\mu_q$
over $\R/\Z$.  Since $\mu_q$ is $1/q$-periodic and $\|\mu_q\|_{\epsilon_2}
\leq e^{\delta_1 q}$ with small $\delta_1$,
$\sup_{x \in \R/\Z} |\mu_q(x)-\lambda| \leq e^{-\delta_2 q}$, for some
$\delta_2=\delta_2(\epsilon_2)$.  Since $|\mu_q(x_0)^{-1}| \geq e^{-\delta_1
q}$ for each $x_0 \in \R$, this implies that $|\lambda|
\geq e^{-\delta_1 q}/2$.  In particular, $\sup_{x \in \R/\Z} |\mu_q(x)-\lambda|<|\lambda|/2$,
so that $d=0$.

Let $\phi^{(k)}=\sum_{j=0}^{k-1} \phi(x+j p/q)$.
By hypothesis,
$\|e^{2 \pi i \phi^{(k)}}\|_{\epsilon_2'} \leq e^{\delta_1 q}$, $1 \leq k
\leq q$ and $\|e^{-2 \pi i \phi^{(q)}}\|_{\epsilon_2'} \leq e^{\delta_1 q}$.
This implies that
\be
\|\phi^{(k)}-k \hat \phi_0\|_{\epsilon_3} \leq C \delta_1 q.
\ee
Indeed, through $\{|\Im z|<\epsilon_2'\}$, it is
obvious that $-\Im \phi^{(k)} \leq \delta_1 q$, $1 \leq k \leq q$ and $\Im
\phi^{(q)} \leq \delta_1 q$, and since for $1 \leq k \leq q-1$ we have
$\phi^{(k)}(z)+\phi^{(q-k)}(z+k p/q)=\phi^{(q)}(z)$, we can conclude that
we have the estimate $|\Im \phi^{(k)}(z)|=2 \delta_1 q$ through the same
band. To estimate the real part, one just uses that
harmonic conjugation in $\{|\Im z|<\epsilon_2'\}$
composed with restriction to $\{|\Im z|<\epsilon_3\}$
is a bounded operator on bounded harmonic functions.

For every $k \in \Z \setminus q \Z$, let $1 \leq j_k \leq q-1$ be such
that $|1-e^{2 \pi i j_k k p/q}| \geq |1-e^{2 \pi i/3}|$.  Then
\be
\frac {|\hat \phi_k|} {|1-e^{2 \pi i k p/q}|}=
\frac {|\hat \phi^{(j_k)}_k|} {|1-e^{2 \pi i j_k k p/q}|} \leq C \delta_1 q
e^{-2 \pi \epsilon_3 |k|}.
\ee
Let
\be
\psi(x)=\sum_{k \in \Z \setminus q \Z}
\frac {\hat \phi_k} {e^{2 \pi i k p/q}-1},
\ee
so that $\|\psi\|_\epsilon \leq C \delta_1 q$.
Let $\theta=\phi^{(q)}/q$, so that
$\|\theta\|_{\epsilon_3}=\|\phi^{(q)}\|_{\epsilon_3}/q \leq C \delta_1$.
Then $\phi(z)=\theta(z)+\psi(z+p/q)-\psi(z)$ (check the Fourier series).
Since
$\|\theta-\hat \theta_0\|_{\epsilon_3} \leq C \delta_1$
and $\theta$ is $1/q$-periodic, we have
$\|\theta-\hat \theta_0\|_\epsilon \leq e^{-\delta q}$ (check the
Fourier series).

The last statement follows automatically from the construction.
\end{pf}

\subsection{Construction of the conjugacy}


Fix $\epsilon<\epsilon_2<\epsilon_1<\epsilon_0$.  Below, $C$ is a large
constant depending on $\epsilon,\epsilon_2,\epsilon_1,\epsilon_0$ that may
increase (finitely many times) along the argument.
Clearly $\ln \|t\|_{\epsilon_0} \leq \ln (2 \|A_q\|_{\epsilon_0}) \leq
\delta_1 q+\ln 2$.
Since $t$ is $1/q$-periodic, it easily follows (by considering the
Fourier series), that $\|t-\hat t_0\|_{\epsilon_1}<e^{-\delta_3 q}$ for
some $\delta_3=\delta_3(\epsilon_0,\epsilon)$.

In particular, $t(x)=\lambda(x)+\lambda(x)^{-1}$ with
$\lambda \in C^\omega_{\epsilon_1}(\R/\Z,\C)$.  Notice also
that $|\lambda(x)|=1$ for $x \in \R$.

Applying Lemma \ref {Pu=0} to $P=A_q-\lambda \id$,
we conclude that there exists $u \in C^\omega_{\epsilon_2}(\R/\Z,\C)$
with $0 \leq -\ln \|u(z)\| \leq C \delta_1 q$, $|\Im z|<\epsilon_2$,
such that $A_q(z) \cdot u(z)=\lambda(z) u(z)$.  Notice that $A(z)
\cdot u(z)$ is a multiple of $u(z+p/q)$ for every $z$, $A(z)=\mu(z)
u(z+p/q)$.  Let $\mu_k$ be as in Lemma \ref {mu}.  We clearly have
\be
\max_{1 \leq k \leq q}
\ln \|\mu_k\|_{\epsilon_2}, \ln \|\mu_k^{-1}\|_{\epsilon_2}
\leq C \delta_1 q+\ln \|A_k\|_{\epsilon_2} \leq C \delta_1 q.
\ee
Let $\psi$ and $\theta$ be given by Lemma \ref {mu},
and let $v=e^{2 \pi i \psi} u$.
Then $A(z) v(z)=e^{2 \pi i \theta(z)} v(z+p/q)$.  Notice that
$-C \delta_1 q \leq \ln \|v(z)\| \leq C \delta_1 q$ through $\{|\Im
z|<\epsilon\}$ and $\|\theta-\hat \theta_0\|_\epsilon \leq e^{-\delta_4 q}$
for some $\delta_4=\delta_4(\epsilon_2,\epsilon)>0$.

Let $\tilde B(z)$ be the matrix with columns $v(z)+\overline
{v(\overline z)}$ and $\frac {1} {i} (v(z)-\overline {v(\overline z)}$.
Then $A(z) \tilde B(z)=\tilde B(z+p/q) R_{\theta(z)}$.
In particular, $b(z+p/q)=b(z)$, where $b=\det \tilde B$.  Since
$\ln \|b\|_{\epsilon_2} \leq C \delta_1 q$, we conclude that
$\|b-\hat b_0\|_\epsilon \leq e^{-\delta_5 q}$, for some
$\delta_5=\delta_5(\epsilon_2,\epsilon)$, provided $\delta_1$ is
sufficiently small.

We claim that $\ln \|b^{-1}\|_\epsilon \leq C \delta_1 q$.  Since $b$ has
exponentially small oscillation, it suffices to show that $\ln b(x_0)^{-1}
\leq C \delta_1 q$ for some $x_0 \in \R$.  If this does not hold, then
there exists
$\kappa \in \C$ with $|\kappa|=1$ such that $-\ln \|v(x_0)-\kappa
\overline {v(x_0)}\| \gg \delta_1 q$.  But since $A_q(x_0)
(v(x_0)-\kappa \overline {v(x_0)})=\lambda(x_0) v(x_0)-\kappa
\lambda(x_0)^{-1} \overline {v(x_0)}$, we conclude that
$-\ln |\lambda(x_0)-\lambda(x_0)^{-1}|-\ln \|v(x_0)\| \gg \delta_1 q$. 
Since $-\ln (2-|t(x_0)|) \leq C \delta_1 q$, this implies that
$-\ln \|v(x_0)\| \gg \delta_1 q$, contradiction.

Up to changing $\theta$ to $-\theta$ (and replacing $v(z)$ by $\overline
{v(\overline z)}$, we may assume that $b(x)>0$ for $x \in \R$.
Now let $B^{-1}=\frac {1} {b^{1/2}} \tilde B$.  Then $\|B\|_\epsilon \leq e^{C
\delta_1 q}$ and $B(x+p/q) A(x)
B(x)^{-1}=R_{\theta(x)}$, as desired.
\qed

\section{Proof of Theorem \ref {q}}

\subsection{Preliminary estimates}

An important input in our estimates is the polynomial bound on solutions of
the Corona problem.  Those can already be found in the original work of
Carleson \cite {C}, but the more precise version given here
has been proved using Wolff's approach.

\begin{thm}[Uchiyama \cite {U}, see Trent \cite {T}]

There exists $C>0$ with the following property.
Let $a_i:\D \to \C$, $1 \leq i \leq k$,
be such that $\delta \leq (\sum_{i=1}^k |a_i(x)|^2)^{1/2}
\leq 1$, $x \in \D$.  Then there
exists $\tilde a_i:\D \to \C$ such that $(\sum_{i=1}^k |\tilde a_i|^2)^{1/2}
\leq C \delta^{-2} (1-\ln \delta)$ and such that
$\sum_{i=1}^k a_i \tilde a_i=1$.

\end{thm}


It is easy to see that Corona estimates for the disk
easily imply corresponding ones for the annulus.

\begin{lemma} \label {K_0}

Let $a_i \in C^\omega_\delta(\R/\Z,\C)$, $1 \leq i \leq k$,
be such that $\delta \leq (\sum_{i=1}^k |a_i(z)|^2)^{1/2} \leq 1$ through
$\{|\Im z|<\delta\}$.  Then there exists $\tilde a_i \in
C^\omega_\delta(\R/\Z,\C)$, $1 \leq i \leq k$,
such that $(\sum_{i=1}^k |\tilde a_i(z)|^2)^{1/2} \leq C
\delta^{-2} (1+\ln \delta)$ (with $C>0$ as in the previous theorem)
and such that $\sum_{i=1}^k a_i \tilde a_i=1$.
Moreover, if all the $a_i$ are real-symmetric, we can choose all the
$\tilde a_i$ real-symmetric.

\end{lemma}

\begin{pf}

By the previous theorem, there exist $a'_i \in C^\omega_\epsilon(\R,\C)$
be such that $\sum_{i=1}^k a_i a'_i=1$ and $(\sum_{i=1}^k
|a'_i(z)|^2)^{1/2} \leq C
\delta^{-2} (1+\ln \delta)$ (we use that the strip $\{|\Im z|<\delta\}$
is conformally equivalent to $\D$).
Let $a^{(j)}_i(z)=\frac {1}
{j} \sum_{n=0}^{j-1} a'_i(z+n)$.  Let $j_n \to \infty$ be a sequence
such that for every $1 \leq i \leq k$, $a^{(j_n)}_i$ converges in the
topology of uniform convergence on compact sets, and let $\tilde a_i$ be the
limits.  Then $\tilde a_i \in C^\omega_\delta(\R/\Z,\C)$.

For the last statement, notice that if the $a_i$ are real-symmetric then
we can substitute each
$\tilde a_i(z)$ by $\frac {1} {2} (\tilde a_i(z)+
\overline {\tilde a_i(\overline z)})$.
\end{pf}




\begin{lemma} \label {detP0}

There exists $C>0$ with the following property.
Consider a function
$P^{(0)}=\left (\bm a^{(0)}&b^{(0)}\\c^{(0)}&d^{(0)} \em \right )$
with coordinates in $C^\omega_{\epsilon_0}(\R/\Z,\C)$.
Assume that $\delta \leq \|P^{(0)}(z)\| \leq 1$ through $\{|\Im
z|<\epsilon_0\}$ and let $\|\det P^{(0)}\|_{\epsilon_0}=\rho
(\delta^2/\ln \delta)^2$.
If $\rho<C^{-1}$ then there exists $P=\left (\bm a&b\\c&d \em \right )$
with coordinates $C^\omega_{\epsilon_0}(\R/\Z,\C)$
such that $\|P^{(0)}-P\| \leq
-C \rho (\delta^2/\ln \delta)$ and $\det P=0$. 
Moreover, if $P^{(0)}$ is real-symmetric then $P$ can be chosen
real-symmetric.

\end{lemma}

\begin{pf}

Let $K_0=\left (\bm a&b\\c&d \em \right )$ be such
that $a^{(0)}d+ad^{(0)}-b^{(0)}c-bc^{(0)}=1$.
Let $P^{(1)}=P^{(0)}-K_0 \det P^{(0)}$.
Then $\det P^{(1)}=(\det P^{(0)})^2 \det K_0$.  Choosing $K_0$ with
minimal $\|K_0\|_{\epsilon_0}$, using Lemma \ref {K_0},
we get $|\det P^{(1)}| \leq C
\rho^2 (\delta^2/\ln \delta)^2$, while $\|P^{(1)}-P^{(0)}\|_{\epsilon_0}<
-C \rho \delta^2/\ln \delta$.  Iterating this
procedure we get a sequence $P^{(n)}$ converging to $P$ as desired.
\end{pf}

\begin{lemma} \label {wreal}

Let $w \in C^\omega_\epsilon(\R/\Z,\C^2)$ be such that
$\delta \leq \|w(z)\| \leq 1$ through $\{|\Im z|<\epsilon\}$.
If $w(x)$ is a multiple of a real vector for $x \in \R/\Z$,
then there exists $\tilde w \in C^\omega_\epsilon(\R,\R^2)$ such that
$\tilde w(x+1)=\pm \tilde w(x)$ for every $x \in \R$,
$w(x)$ is a multiple of $\tilde w(x)$ for every $x \in \R$,
and $C^{-1} \delta^{3/2} \leq \|\tilde w(x)\| \leq C \delta^{-1/2}$
for $|\Im x|<\epsilon$.

\end{lemma}

\begin{pf}

Let $w=(a,b)$, and
let $\tilde a(z)=\overline {a(\overline z)}$, $\tilde b(z)=\overline
{b(\overline z)}$.  Let $\phi=a/\tilde a=b/\tilde b$.  Then $2^{-1/2}
\delta \leq \phi(z) \leq 2^{1/2} \delta^{-1}$ through
$\{|\Im z|<\epsilon\}$.  Let $\tilde w(x)=\phi^{-1/2} w(x)$.
\end{pf}

\begin{rem} \label {real}

If in the statement of Lemma \ref {Pu=0} we further assume that $P$ is
real-symmetric, we can then obtain, using Lemma \ref {wreal},
a real-symmetric solution of $Pu=0$ satisfying the required bounds (with
adjusted constants), which is not necessarily
$1$-periodic, but satisfies $u(x)=\pm u(x+1)$.  It is not possible in
general to get a $1$-periodic solution, as exemplified by $P(x)=\left (\bm
-\sin \pi x \cos \pi x&-\sin^2 \pi x \\ \cos^2 \pi x & \sin \pi x \cos \pi x
\em \right )$.

\end{rem}

\begin{lemma} \label {delta12}

For every $0<\epsilon_1<\epsilon_0$, there exists $C_0>0$, $\delta_0>0$
such that for every $0<\delta<\delta_0$, if $q$ is sufficiently large and
$\phi \in C^\omega_{\epsilon_0}(\R/\Z,\C)$ is such that
$\|\phi\|_{\epsilon_0} \leq
e^{\delta q}$, while there exists $z$ with $|\Im z|<\epsilon_1$ such that
$|\phi(z)|,...,|\phi(z+(q-1)/q)|<e^{-C \delta q}$ with $C>C_0$, then
$\|\phi\|_{\epsilon_1} \leq \max \{e^{-C_0^{-1} C \delta q},e^{-\delta_0
q}\}$.

\end{lemma}

\begin{pf}

Let $\tilde \phi(x)=\sum_{k=-[q/2]}^{q-1-[q/2]} \hat \phi_k e^{2 \pi i k
x}$.  Then $\|\phi-\tilde \phi\|_{\epsilon_1} \leq e^{-\delta_1 q}$ with
$\delta_1=\delta_1(\epsilon_0,\epsilon)>0$, provided $\delta_0$ is
sufficiently small.  By Lagrange interpolation,
\be
\sup_{\Im x=\Im z} \|\tilde \phi(x)\| \leq \sum_{k=0}^{q-1}
|\phi(z+k/q)| \leq q e^{-C \delta q}.
\ee
(Indeed $\psi(x)=
e^{2 \pi i [q/2] x} \tilde \phi(x+z)$ satisfies
$\psi(x)=\sum_{k=0}^{q-1} \psi(k/q) c_q(x-k/q)$
where $c_q(x)=\frac {1} {q} \sum_{k=0}^{q-1} e^{2 \pi i k x/q}$, so that
$\sup_{x \in \R/\Z} |c_q(x)|=1$.)
The result follows by convexity.
\end{pf}

\subsection{Construction of the conjugacy}

Fix $\epsilon<\epsilon'<\epsilon_0$.
As in the proof of Theorem \ref {qsimpler}, the hypothesis
implies that
\be \label {d3}
\|t(x)-\hat t_0\|_{\epsilon'} \leq e^{-\delta_3 q}
\ee
for some $\delta_3=\delta_3(\epsilon',\epsilon_0)>0$.

We will consider $2$ distinct cases:
\begin{enumerate}
\item $|2-|\hat t_0|| \geq e^{-C_0^2 \delta_1 q}$,
\item $|2-|\hat t_0||<e^{-C_0^2 \delta_1 q}$.
\end{enumerate}
Here $C_0=C_0(\epsilon_0,\epsilon)$ will be some appropriately large
constant.



\comm{
Let us assume by contradiction that there exists a sequence $p_n/q_n$ and
$A^{(n)}$ such that the result does not hold.  Up to considering
subsequences, we may assume to eventually fall into one of the following
cases:
\begin{enumerate}
\item $|\hat t_0|<2$ and $2-|\hat t_0| \geq e^{-o(q)}$,
\item $|\hat t_0|>2$ and $|\hat t_0|-2 \geq e^{-o(q)}$,
\item There exists $|2-|\hat t_0||$ is exponentially small but $\|A_q(z)-\id\|_\epsilon
\geq e^{-o(q)}$.
}

\subsubsection{Case (1)}

If $|\hat t_0|<2$, the result is contained in Theorem \ref {qsimpler}.  We
will assume below that $|\hat t_0|>2$.  In this case (\ref {d3}) implies
$t(z)=\lambda(z)+\lambda(z)^{-1}$ with $\lambda \in
C^\omega_{\epsilon'}(\R/\Z,\R)$.  The argument below is a simple adaptation
of that of Theorem \ref {qsimpler}.  Fix
$\epsilon<\epsilon''<\epsilon'''<\epsilon_1$.

\comm{
Applying Lemma \ref {Pu=0} to $P=A_q-\lambda \id$,
we conclude that there exists $u \in C^\omega_{\epsilon''}(\R/\Z,\C^2)$ with
$e^{-o(q)} \leq \|u(x)\| \leq 1$, $|\Im x|<\epsilon''$,
such that $A_q \cdot u=\lambda u$.  Notice that $A(x)
\cdot u(x)$ is a multiple of $u(x+p/q)$ for every $x$, $A(x)=\mu(x)
u(x+p/q)$.  Let $\mu_k$ be as in Lemma \ref {}.  We clearly have
$\|\mu_k\|_{\epsilon''}, \|\mu_k^{-1}\|_{\epsilon''}
\leq e^{o(q)} \|A_k\|_{\epsilon''} \leq e^{o(q)}$,
$1 \leq k \leq q$.  Let $\psi$ and $\theta$ be given by Lemma \ref {}(with
$\epsilon_0=\epsilon''$ and $\epsilon=\epsilon'''$), and let
$v(x)=e^{2 \pi i \psi(x)} u(x)$.  Then $A(x) v(x)=e^{2 \pi i \theta(x)}
v(x+p/q)$.
}


Applying Lemma \ref {Pu=0} and Remark \ref {real} to $P=A_q-\lambda \id$, we
obtain $u \in C^\omega_{\epsilon''}(\R,\R^2)$ with $e^{-C \delta_1 q} \leq
\|u(x)\| \leq 1$, $|\Im x|<\epsilon''$ such that
$u(z+1)=\pm u(z)$ and $A_q u=\lambda u$.
Notice that
$A(z) \cdot u(z)$ is a multiple of $u(z+p/q)$ for every $z$,
$A(z) u(z)=\mu(z) u(z+p/q)$.  Notice that $\mu$ is $1$-periodic.
Let $\mu_k$ be as in Lemma \ref {mu}.  We clearly have
$\|\mu_k\|_{\epsilon''}, \|\mu_k^{-1}\|_{\epsilon''}
\leq e^{C \delta_1 q} \|A_k\|_{\epsilon''} \leq e^{C \delta_1 q}$,
$1 \leq k \leq q$.  Let $\psi$ and $\theta$ be given by Lemma \ref {mu}
(with $\epsilon_0=\epsilon''$ and $\epsilon=\epsilon'''$), and let
$v=e^{2 \pi i \psi} u$.  Then $v$ is real-symmetric
and $A(z) v(z)=\gamma(z)
v(z+p/q)$, where $\gamma(z)$ is $1/q$-periodic and real-symmetric.  Notice
that $\gamma$ has the same sign as $\mu$ and $\gamma^q=\lambda$.

An analogous argument yields a solution $v' \in C^\omega_{\epsilon'''}(\R,\R^2)$
such that $v'(z+1)=\pm v(z)$, $A_q v'=\lambda^{-1} v'$ and $A(z)
v'(z)=\gamma(z)^{-1} v(z+p/q)$.
Notice that since $\lambda \neq \lambda^{-1}$, $v$ is not colinear with
$v'$, so the determinant of the matrix with columns $v(x)$ and $v'(x)$ does not
change sign for $x \in \R$.  Thus
$v$ and $v'$ are both $1$-periodic or $1$-antiperiodic. 
Take $\tilde B$ as the matrix with columns $v$ and $v'$.  Since $A(z) \tilde
B(z)=\tilde B(z+p/q) \bp\gamma(z)&0\\0&\gamma(z)^{-1} \ep$, $b=\det \tilde B$ is
$1/q$-periodic, so that $\|b-\hat b_0\|_\epsilon \leq e^{-\delta_2 q}$
(where $\delta_2>0$ is independent of $\delta_1$ small).
For fixed $x_0 \in \R$, since $A_q(x_0)$ is an $e^{\delta_1 q}$
bounded matrix whose eigenvalues
$\lambda(x_0)$ and $\lambda^{-1}(x_0)$ are $e^{-C \delta_1 q}$ apart,
the angle between the eigenvectors $v(x_0)$ and $v'(x_0)$ is at least $e^{-C
\delta_1 q}$. 
Thus $|b(x_0)| \geq e^{-C \delta_1 q}$, and hence
$|\hat b_0| \geq e^{-C \delta_1 q}$. 
The result follows by taking $B^{-1}$ as the matrix with columns $v$ and $v'/b$.

\subsubsection{Case (2)}


We will use, in two distinct situations, the following estimate.

\begin{lemma} \label {WR}

For every $\epsilon<\epsilon_1<\epsilon_0$,
there exists $C_5>1$ such that for every $C_3>1$ sufficiently large and
every $C_4>1$, for every $\delta_1>0$ is sufficiently small, if
$p/q \in \Q$ with $q$ sufficiently large, and
$A \in C^\omega_{\epsilon_0}(\R/\Z,\SL(2,\R))$ satisfies (\ref {cond}), then
the following property holds.  If there exists
$W \in C^\omega_{\epsilon_1}(\R/\Z,\SL(2,\R))$ and $R \in \SO(2,\R)$
satisfying
\be \label {R}
\|W(z+p/q) A(z)-R W(z)\|_{\epsilon_1} \leq e^{-C_3 C_4 \delta_1 q},
\ee
while
\be \label {Wz}
\|W(z)\| \geq e^{-C_4 \delta_1 q}, \quad |\Im z|<\epsilon_1,
\ee
and
\be
\|\det W\|_{\epsilon_1} \leq e^{-C_3 C_4 \delta_1 q},
\ee
then there exists $B \in
C^\omega_\epsilon(\R/\Z,\PSL(2,\R))$ and a constant diagonal matrix $D \in
\SL(2,\R)$ such that $\|B\|_\epsilon \leq e^{C_5 C_4 \delta_1 q}$ and
$\|B(z+p/q) A(z) B(z)^{-1}-D\|_\epsilon \leq e^{-C_4 \delta_1 q}$.

\end{lemma}

\begin{pf}

Fix $\epsilon<\epsilon_2<\epsilon_1$.
Apply Lemma \ref {detP0} to $P^{(0)}=W$ to obtain $P$ with $\det P=0$
such that
$\|P-W\|_{\epsilon_1} \leq e^{C^{-1} C_3 C_4 \delta_1 q}$
Using Lemma \ref {Pu=0} together with Remark \ref {real}
to get $u \in C^\omega_{\epsilon_2}(\R,\R^2)$
such that $u(z+1)=u(z)$ or $u(z+1)=-u(z)$ such that $P u=0$, and satisfying
$e^{-C C_4 \delta_1 q} \leq \|u(z)\| \leq 1$ through $\{|\Im
z|<\epsilon_2\}$.
Then
\be
\|W(z) \cdot u(z)\|_{\epsilon_2} \leq e^{-C^{-1} C_3 C_4 \delta_1 q}.
\ee
Using (\ref {R}) we get
\be
\|W(z+p/q) \cdot A(z) u(z)\|_{\epsilon_2} \leq e^{-C^{-1}
C_3 C_4 \delta_1 q},
\ee
Using (\ref {Wz}), it follows that
$A(z) \cdot u(z)$ is $e^{-C^{-1} C_3 C_4 \delta_1 q}$ close to a multiple
of $u(z+p/q)$, $|\Im z|<\epsilon_2$.
Using Lemma \ref {K_0}, define
$\tilde B \in C^\omega(\R/\Z,\PSL(2,\R))$ with first column $u$.  Then
\be
\tilde B(z+p/q)^{-1} A(z) B(z)=
\left (\bm \mu(z)&\tilde s_2(z)\\\tilde s_3(z)&\mu(z)^{-1}+\tilde s_4(z) \em
\right ),
\ee
with $\mu$ real-symmetric, and
$\|\tilde s_3\|_{\epsilon_2},\|\tilde s_4\|_{\epsilon_2}<e^{-C^{-1} C_3 C_4 \delta_1
q}$.  As in the proof of case (2) above,
we can apply Lemma \ref {mu} to obtain $\psi$ and
$\theta$ such that $B'=\bp e^{-\psi}&0\\0&e^{\psi} \ep \tilde B$ satisfies
\be
B'(z+p/q) A(z) B'(z)^{-1}=\left (\bm e^{2 \pi i \theta(z)}
&s_2'(z)\\s_3'(z)&e^{-2 \pi i \theta(z)}+s_4'(z) \em
\right ),
\ee
with $\|s_3'\|_\epsilon,\|s_4'\|_\epsilon \leq e^{-C^{-1} C_3 C_4 \delta_1
q}$.  We also have the bound $\|B'\|_\epsilon \leq e^{C C_4 \delta_1 q}$,
and hence $\|s_2'\|_\epsilon \leq e^{C C_4 \delta_1 q}$.  Since
$\|\theta-\hat \theta_0\|_\epsilon \leq e^{-\delta q}$ with
$\delta=\delta(\epsilon,\epsilon_2)$, the result follows
with $B^{-1}=\bp
d&0\\0&d^{-1} \ep B'$, $d=e^{10 C_4 \delta_1 q} (1+\|s_2'\|_\epsilon)$.
\end{pf}

One key consideration when $\hat t_0$ is close to $\pm 2$ is whether $\pm A_q$ is
close to the identity or not.  Fix $\epsilon<\epsilon_1<\epsilon'$.
Notice that if
$\|A_q \mp \id\|_{\epsilon_1} \geq e^{-C_0 \delta_1 q}$ then
\be \label {Aq-id}
\|A_q(z) \mp \id\| \geq e^{-C_1 C_0 \delta_1 q}, \quad |\Im z|<\epsilon_1,
\ee
for an appropriately large constant $C_1$, which does not depend on the
choice of $C_0$.
Indeed, if this was not the case then there would exists $z$ with $|\Im
z|<\epsilon_1$ such that
$\|A_q(z+k p/q) \mp \id\| \leq \|A_k(z)\|^2 \|A_q(z) \mp \id\| \leq
e^{-C_2 C_0 \delta_1 q}$ for $0 \leq k \leq q-1$ with $C_2$ large.
Applying Lemma \ref {delta12} to the coefficients of $A_q \mp \id$, that are
bounded by $e^{\delta_1 q}+1$ through $\{|\Im z|<\epsilon_0\}$, leads to a
contradiction.

We will assume below that $C_0$ is chosen much bigger than $C_1$.
Then, under the assumption that (\ref {Aq-id}) holds, the result follows
from Lemma \ref {WR}, with $W=A_q \mp \id$.

Assume not that (\ref {Aq-id}) does not hold, so that, as explained above,
we must have
\be \label {aq}
\|A_q \mp \id\|_{\epsilon_1} \leq e^{-C_0 \delta_1 q}.
\ee

Let us consider a large coefficient
of the discrete Fourier transform of the essentially periodic sequence
$\{R_{l s/2 q} A_s\}_{s=0}^{q-1}$, where $l=0$ if $A_q$ is close to $\id$ and
$l=1$ is $A_q$ is close to $-\id$.
More precisely, take $W_k=\sum_{s=0}^{q-1} R_{k s/q} R_{l s/2 q} A_s$, $0 \leq k
\leq q-1$.  Then
\be \label {Wk}
W_k(z+p/q) A(z)=R_{-(2 k+l)/2 q}(W_k(z) \pm A_q(z)-\id),
\ee
so that by (\ref {aq}),
\be \label {3}
\|W_k(z+p/q) A(z)-R_{(2 k+l)/2 q} W_k(z)\|_{\epsilon_1} \leq e^{-C^{-1}
C_0 \delta_1 q}
\ee
(here and below, we use $C$ for quantities that do not become larger if
$C_0$ is taken large, so that we can always assume that $C_0>C$).
Clearly, for every $x \in \R/\Z$ and any unit vector $y \in \R^2$,
\be
\sum_{k=0}^{q-1} \|W_k(x) \cdot y\|^2=q \sum_{s=0}^{q-1} \|A_s(x) \cdot y\|^2
\ee
(Parseval identity).
The average of the right hand side over the circle of unit vectors is
\be
q \sum_{s=0}^{q-1} \frac {\|A_s(x)\|^2+\|A_s(x)\|^{-2}} {2} \geq q^2,
\ee
so for every $x \in \R/\Z$, there exists one unit vector $y$ such that
$\sum_{k=0}^{q-1} \|W_k(x) \cdot y\|^2 \geq q^2$.
Fix $x_0 \in \R$ and let $W=W_{k_0}$ where $k_0$ is such that
such that $\|W_{k_0}(x_0)\|^2$ is maximal.  Then $\|W(x_0)\|^2 \geq q$.

We claim that for fixed $\epsilon<\epsilon_1'<\epsilon_1$,
\be \label {w(z)}
-\ln \|W(z)\| \leq C \delta_1 q, \quad |\Im z|<\epsilon_1'.
\ee
Indeed, if $\|W(z)\| \leq e^{-C \delta_1 q}$ with large $C$, then
$\max_{0 \leq j \leq q-1} \|W(z+j p/q)\| \leq e^{-C \delta_1 q}$ with large
$C$, and Lemma \ref {delta12} implies that $\|W\|_{\epsilon_1'} \leq e^{-C
\delta_1 q}$ with large $C$.
This contradicts $\|W(x_0)\|^2 \geq q$, so that (\ref {w(z)}) holds.

Let $w=\det W$.  Then, by (\ref {3}),
$\|w(z+p/q)-w(z)\|_{\epsilon_1} \leq e^{-C^{-1} C_0
\delta_1 q}$, so $w$ is $e^{-C^{-1} C_0 \delta_1 q}$ to $\frac {1} {q}
\sum_{k=0}^{q-1} w(x+k p/q)$, which is $1/q$-periodic and bounded by
$q e^{2 \delta_1 q}$ over
$\{|\Im z|<\epsilon_1\}$.  By convexity, we get that
$\|w-\hat w_0\|_{\epsilon_1'} \leq e^{-C^{-1} C_0 \delta_1 q}$.

Assume that $-\ln |\hat w_0| \leq C_0^{1/2} \delta_1 q$.
If $w(x)>0$ for $x \in
\R$, then take $B=\frac {1} {w^{1/2}} W$.  If $w(x)<0$ for $x \in
\R$, then take $B=\bp 1&0\\0&-1 \ep
\frac {1} {(-w)^{1/2}} W$.  The result follows.

Assume that $-\ln |\hat w_0| \geq C_0^{1/2} \delta_1 q$.
The result follows by applying Lemma \ref {WR} (with $\epsilon_1'$ instead of
$\epsilon_1$).
\qed

\begin{rem}

The analysis above can be refined further to yield considerably more precise
estimates.

\end{rem}

\end{document}